\documentclass[11pt, letterpaper]{article}
\usepackage[utf8]{inputenc}
\usepackage[margin=1.5cm]{geometry}
\usepackage{titlesec}
\usepackage{tabu}
\usepackage{enumitem}
\usepackage{amssymb}
\usepackage{graphicx}
\usepackage{xcolor}
\newlist{selectlist}{itemize}{2}
\setlist[selectlist]{label=$\square$,leftmargin=*,noitemsep,topsep=0pt}
\usepackage{ragged2e}
\usepackage{amsmath} 
\usepackage{lmodern}
\usepackage{float}
\usepackage{hyperref}
\hypersetup{
    colorlinks=true,
    linkcolor=blue,
    filecolor=magenta,      
    urlcolor=blue,
}
 
\titleformat{\section}[block]{\hspace{1em}\bfseries}{\thesection.}{0.5em}{} 
\titleformat{\subsection}[block]{\hspace{1em}}{\thesubsection}{0.5em}{}

\begin{document}

\begin{center}
\noindent
\textbf{TumorPred: A Computational Framework Implemented via an R/Shiny Web Application for Parameter Estimation and Sensitivity Analysis in Compartmental Brain Modeling}
\end{center}
\vskip0.5cm

\begin{center}
\noindent
\textbf{Charuka D. Wickramasinghe$^{a, *}$, Nelum S. S. M. Hapuhinna$^{b}$\\ 
}
\end{center}
\vskip0.1cm

\noindent \textbf{\textit{\small $^{a}$Karmanos Cancer Institute, School of Medicine, Wayne State University, Detroit, MI, 48202, USA}}
\textbf{\textit{\small $^{b}$Department of Mathematics and Statistics, Northern Kentucky University, KY, 41099, USA}}\\
\textbf{\textit{\small $^{*}$Corresponding author. gi6036@wayne.edu}}\\

\noindent
\textbf{Abstract}\\
It is difficult or infeasible to directly measure how much of a drug actually enters the human brain and a brain tumor, how long it remains there, and to estimate drug-specific or patient-specific parameters, as well as how changes in these parameters influence model outputs and pharmacokinetic characteristics. Compartmental modeling offers a powerful mathematical approach to describe drug distribution and elimination in the body using systems of differential equations. This study introduces TumorPred, an R/Shiny-based web application designed for model simulation, sensitivity analysis, and pharmacokinetic parameter calculation in a permeability-limited four-compartment brain model. The model closely mimics human brain functionality for drug delivery and aims to predict the pharmacokinetics of drugs in the brain blood, brain mass, and cranial and spinal cerebrospinal fluid (CSF) of the human brain.  The app provides real-time output updates in response to input modifications and allows users to visualize and download simulated plots and data tables. The computational accuracy of TumorPred is validated against results from the Simcyp Simulator (Certara Inc.). TumorPred is freely accessible and serves as an invaluable computational tool and data-driven resource for advancing drug development and optimizing treatment strategies for more effective brain cancer therapy.
\vskip0.5cm

\vskip0.5cm
\newpage
\noindent
\textbf{Code metadata}\\

\noindent
\begin{tabular}{|l|p{6.5cm}|p{9.8cm}|}
\hline
\textbf{Nr.} & \textbf{Code metadata description} & \textbf{} \\
\hline
C1 & Current code version & \textit{1.0} \\
\hline
C2 & Permanent link to code/repository used for this code version & \textit{\underline{\url{https://github.com/CharukaWick/TumorPred_Version_1.0}}} \\
\hline
C3  & Permanent link to Reproducible Capsule & \textit{\underline{\url{https://github.com/CharukaWick/TumorPred_Version_1.0}}} \\
\hline
C4 & Legal Code License   & \textit{ MIT License } \\
\hline
C5 & Code versioning system used & \textit{git} \\
\hline
C6 & Software code languages, tools, and services used & \textit{R}\\
\hline
C7 & Compilation requirements, operating environments \& dependencies & \textit{None}\\
\hline
C8 & If available Link to developer documentation/manual & \textit{N/A}\\
\hline
C9 & Support email for questions & gi6036@wayne.edu\\
\hline
\end{tabular}\\
\vskip0.5cm
\noindent
\section{Introduction}
\label{sec1}

Quantitative knowledge  of drug  absorption, distribution, metabolism, and excretion (ADME) of many new and existing drugs in the human central nerve system (CNS)  is essential for the development of new drugs and optimal use of existing drugs for brain cancer. However, direct measurement of drug penetration in the human brain and brain tumors is difficult or infeasible given the challenge of sampling and limitation of currently available imaging and analytical technologies. Physiologically based pharmacokinetic (PBPK) modeling offers a mechanism based computational approach for quantitatively predicting the ADME of drugs, as it can incorporate both biological system-specific and drug-specific data into a pharmacokinetic model \cite{Li2,nel,pinn}. These PBPK models consist of complex systems of coupled ordinary differential equations. As a result, there is growing interest in developing user-friendly computational tools to solve such models\cite{Li3,ModVizPop,gpk}. These tools enable multidisciplinary teams in the pharmaceutical industry, including researchers and clinicians without programming expertise, to effectively support decision-making in drug development.

To this end, we developed a computational framework (TumorPred) via an R/Shiny web application \cite{shiny}. TumorPred offers three main functionalities. \textbf{Module 1} for model simulations helps direct calculation of drug concentration profiles and estimation of pharmacokinetic parameters. \textbf{Module 2} for sensitivity analysis is designed to understand which parameters are the most influential in predicting drug concentration profiles in the CNS.  \textbf{Module 3} for parameter estimation helps in estimating drug-specific and patient-specific parameters which is crucial for modeling drug behavior in the body and ensuring the effectiveness and safety of medications, while inaccurate estimates can weaken pharmacokinetic predictions, slowing drug development and clinical practices. 

We incorporated a 4-compartment permeability-limited brain model introduced in \cite{Jamei1} into the TumorPred app, as this model is considered one of the gold standards for predicting drug penetration, allowing a wide range of users to benefit from the TumorPred. As supporting evidence, the same brain model is implemented in the Simcyp Simulator\cite{Jamei}, a widely used commercial software for drug discovery. The novelty of TumorPred lies in its ability to estimate both drug-specific and system-specific model parameters using a parameter estimation algorithm, rather than relying on values from the literature or in vitro–in vivo studies. The app is programmed in R and utilizes Shiny as the web application framework, taking advantages of key R packages \cite{ggplot,desolve,tidy}. 

\section{Software description}

TumorePred is an R/Shiny web applications developed using RStudio, one of the most widely adopted integrated development environments (IDEs) for R programming. The application code is organized in a single file containing three main components: the user interface (\texttt{UI}), the \texttt{server}, and the \texttt{shinyApp} function, as illustrated in Figure~\ref{inputt}(d). The application source code has been publicly available via Git.  

\subsection{Software requirements}
In the pharmaceutical industry, it is important to quantitatively predict how drugs reach and behave in different regions of the human brain using realistic biological structures and physiological conditions. Secondly, due to experimental limitations, accurate estimation of drug-related parameters is essential. Thirdly, sensitivity analysis is important because it helps scientists identify which parameters strongly influence the model’s results. TumorPred facilitates all three of these requirements and assists scientists in the pharmaceutical industry, as well as researchers in oncology and pharmacology, in predicting drug penetration, performing parameter estimation, and conducting sensitivity analysis. These simulation results further enable drug developers to estimate key pharmacokinetic parameters which directly influence dose selection, safety evaluation, efficacy assessment, and regulatory decision-making.

\subsection{Software development process}

The \texttt{TumorPred} application was developed using an iterative and modular prototyping approach. 
The development process consisted of the following stages to ensure correctness, flexibility, and usability:
\vspace{0.3cm}

\noindent \textbf{Stage 1. Graphical User Interface (GUI) Design}:
The graphical user interface (GUI) was first developed to define and organize all user interactions with the computational framework. This included data file upload functionality, parameter specification fields, simulation control options, and output visualization panels. This stage established a clear interaction structure between the user and the underlying computational engine.

\begin{figure}[ht]
\centering
\includegraphics[width=18cm]{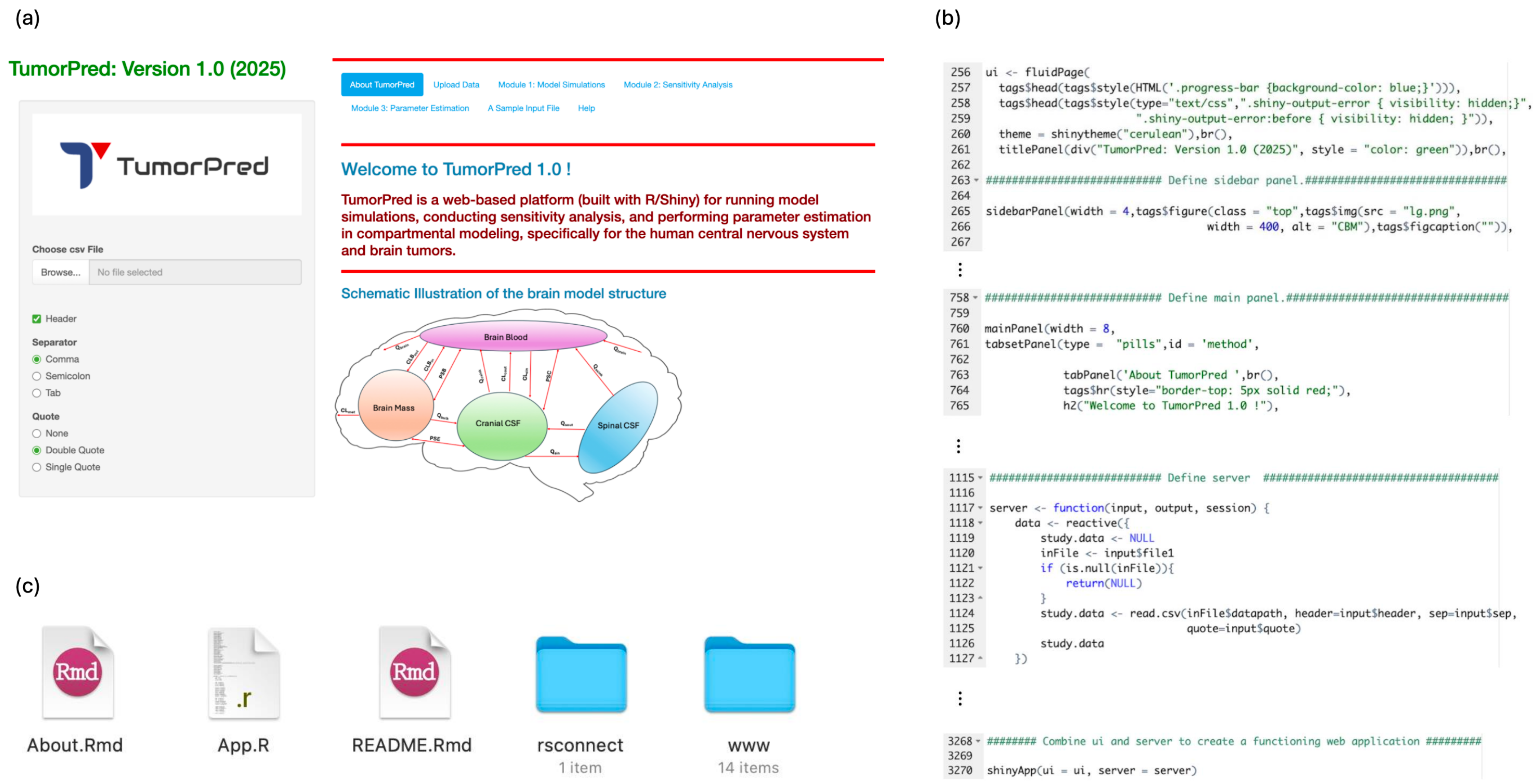}
\begin{tiny}
\caption{The GUI and the code structure (including user-interface and server) of the TumorPred software.}\label{inputt}
\end{tiny}
\end{figure}

The GUI is organized into two primary sections: a sidebar panel (left) for input controls and parameter configuration, and a main panel (right) for simulation outputs and visualizations, as illustrated in Figure~\ref{inputt}(a). Numerical solutions of the system of ordinary differential equations (ODEs) are computed using appropriate ODE solvers to generate simulation results. Detailed parameter descriptions and reference parameter values are provided in the Supplementary Materials.

\vspace{0.3cm}

\noindent \textbf{Stage 2. Server-Side Implementation}:
Inputs provided by the user through the GUI are transferred to the server, where they are processed. The server logic was then developed in a modular manner:
\begin{itemize}
    \item \textbf{Module 1: Implementation of the Brain Model}: The physiologically based brain model governed by a system of four ordinary differential equations (ODEs), which predicts drug penetration into different regions of the human brain, is presented in equations (1)-(4).

{\footnotesize    
\noindent \textbf{\textit{Brain blood compartment:}} 
\begin{equation} \label{eqn1}
\begin{split}
V_{bb} \frac{dC_{bb}}{dt} & = Q_{brain}(C_{art}-C_{bb}) + PSB(\lambda_{bm}fu_{bm}C_{bm}-\lambda_{bb}fu_{bb}C_{bb})  +  CLB_{in}fu_{bb}C_{bb} + CLB_{out}fu_{bm}C_{bm} 
\\ & + PSC(\lambda_{ccsf}fu_{ccsf}C_{ccsf}-\lambda_{bb}fu_{bb}C_{bb}) - CLC_{in}fu_{bb}C_{bb} 
+ CLC_{out}fu_{ccsf}C_{ccsf} + Q_{csink}C_{ccsf} 
+ Q_{ssink}C_{scsf}
\end{split}
\end{equation}
\noindent \textbf {\textit{Brain mass compartment:}}
\begin{equation} \label{eqn2}
\begin{split}
\begin{aligned}
 V_{bm} \frac{dC_{bm}}{dt} &= PSB(\lambda_{bb}fu_{bb}C_{bb}-\lambda_{bm}fu_{bm}C_{bm})  +CLB_{in}fu_{bb}C_{bb}  - CLB_{out}fu_{bm}C_{bm} - Q_{bulk}C_{bm}  - CL_{met}C_{bm} \\
&  + PSE(\lambda_{ccsf}fu_{ccsf}C_{ccsf}-\lambda_{bm}fu_{bm}C_{bm}) 
\end{aligned}
\end{split}
\end{equation}
\noindent \textbf {\textit{Cranial CSF compartment:}}
\begin{equation} \label{eqn3}
\begin{split}
\begin{aligned}
 \qquad V_{ccsf} \frac{dC_{ccsf}}{dt} 
 &= PSC(\lambda_{bb}fu_{bb}C_{bb}-\lambda_{ccsf}fu_{ccsf}C_{ccsf}) 
+CLC_{in}fu_{bb}C_{bb}    - CLC_{out}fu_{ccsf}C_{ccsf}  - Q_{sout}C_{scsf} 
 \\&  + PSE(\lambda_{bm}fu_{bm}C_{bm}-\lambda_{ccsf}fu_{ccsf}C_{ccsf})  - Q_{sin}C_{ccsf} 
 - Q_{csink}C_{ccsf}
\end{aligned}
\end{split}
\end{equation}
\noindent \textbf {\textit{Spinal CSF compartment:}}
\begin{equation} \label{eqn4}
\begin{split}
\begin{aligned}
\qquad V_{scsf} \frac{dC_{scsf}}{dt} &= Q_{sin}C_{ccsf} - Q_{sout}C_{scsf} - Q_{ssink}C_{scsf} \hspace{8cm}
\end{aligned}
\end{split}
\end{equation}
}
    
    \item \textbf{Module 2. Sensitivity Analysis}: Based on the computed ODE solutions, a sensitivity analysis module was developed to evaluate the influence of model parameters on system outputs.
    \item \textbf{Module 3. Parameter Estimation}: A parameter estimation module was incorporated to allow users to estimate unknown model parameters using uploaded experimental or clinical plasma data. 
\end{itemize}

\vspace{0.3cm}

\noindent \textbf{Stage 3. Integration and Output Rendering}: 
The UI and server components were integrated to dynamically render outputs in the main panel of the interface. The command \texttt{shinyApp(ui = ui, server = server)} links the UI and server into a single interactive web application as shown in Figure~\ref{inputt}(b). The server then generates outputs such as figures, tables, or other graphical elements which are displayed in the positions specified by the UI. Furthermore, a \texttt{www} directory (Figure~\ref{inputt}(c)) can be created in the root folder of a Shiny application to store external resources (e.g., scripts, images, excel files, or videos) that can be directly referenced within the app.  The description of the parameters and the reference values can be found in the supplementary materials. Finally, the \texttt{TumorPred} application has been deployed on \texttt{shinyapps.io}, making it accessible to a broader scientific and user community.  

\subsection{Input file preparation and model simulations }
Drug penetration into and distribution within the four CNS compartments are driven by the plasma drug concentration–time profile.  The application allows users to upload structured input files containing time-series plasma data and parameter values. These data are processed within the server to perform simulations, sensitivity analysis, and parameter estimation in a unified workflow. We assume user has a set of state vector observations available for parameter estimation. Therefore, the user input file (.csv) must contain the following data columns: time, plasma, observed data for each compartment, and optionally, a set of parameter values (if the user wishes to perform a sensitivity analysis independently of parameter estimation). A sample input file can be accessed through the GUI of TumorPred; it contains data gathered from the Simcyp Simulator following the oral administration of a single 10 mg dose of abemaciclib drug. Note that the user must retain the original column names in the input file, as the server relies on them to read and process the data upon submission. Upon submitting the data, Module 1 allows users to perform model simulations. The results include an estimated drug concentration data table, concentration-time plots, and key pharmacokinetic parameters (Cmax, Tmax, and AUC), all of which can be viewed and downloaded. For demonstration purposes, a sample input file is available through the GUI; this sample can be used in place of the user's own data.

\subsection{Sensitivity analysis}

Sensitivity analysis is a crucial method for identifying the parameters that most significantly influence the outcomes of a model. In the context of this study, sensitivity analysis can be approached in two ways. First, following a rigorous parameter estimation process, users can update the input file with the estimated parameters to proceed with the sensitivity analysis. Alternatively, if the user already has all the parameter values for a specific drug under investigation, the sensitivity analysis can be conducted directly without the need for parameter estimation. 

\begin{figure}[ht]
\centering
\includegraphics[width=18cm]{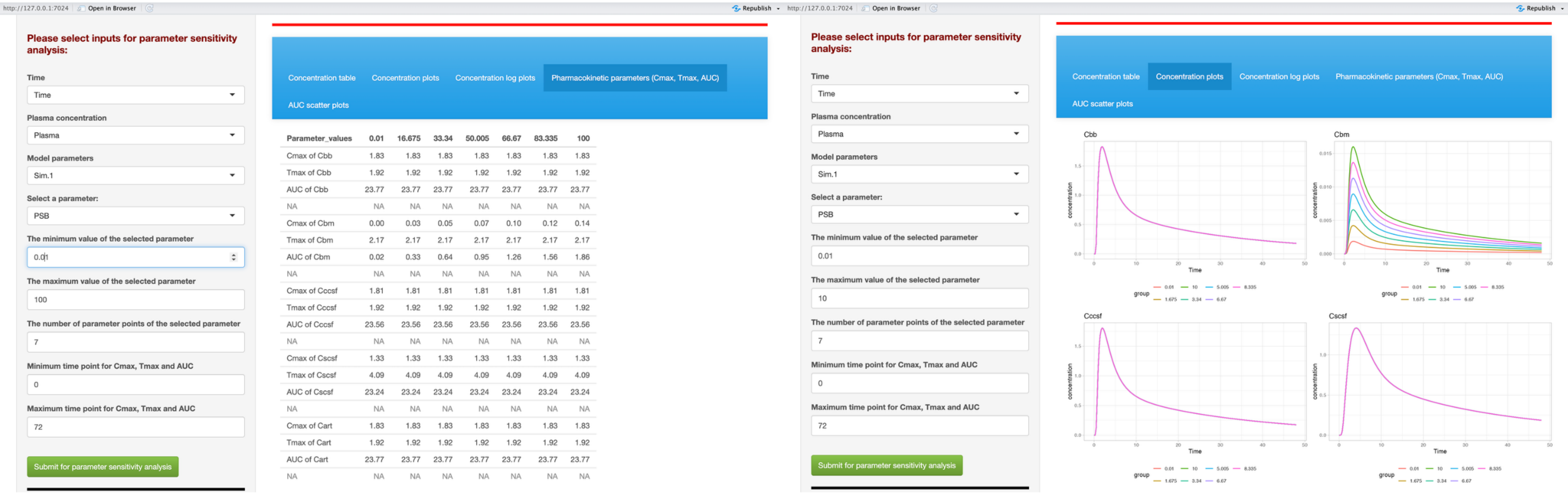}
\begin{tiny}
\caption{Sensitivity analysis results for the parameter PSB, with all other parameters held constant.}\label{inputp}
\end{tiny}
\end{figure}

For the example illustrated in Figure (\ref{inputp}), we investigate the sensitivity of the parameter PSB (passive permeability–surface area product of the blood–brain barrier). To this end, seven values of PSB are generated over the range 0.01 to 100, as shown in the left panel of Figure (\ref{inputp}). For each value, key pharmacokinetic parameters (AUC, Cmax, and Tmax) are computed for the four compartments of the human brain. During this analysis, all other model parameters are held constant to isolate the effect of PSB. In addition, concentration–time profiles are predicted for each of the four brain compartments corresponding to the seven PSB values. The results (Figure \ref{inputp}, right panel) indicate that the brain mass compartment (C\textsubscript{bm}) is highly sensitive to variations in PSB, whereas the remaining compartments (C\textsubscript{bb}, C\textsubscript{ccsf}, and C\textsubscript{scsf}) show minimal or no observable sensitivity to this parameter.

\subsection{Parameter estimation}

To find the best-fit parameter values for a given drug and patient population, we define a loss function ($\text{loss} = \min_{p} \sum_{i=1}^{n} \sum_{j=1}^{k}\frac{(\hat{y}_{j}(t_{i}) - y_{j}(t_{i}, p))^2}{2\sigma^2_{j}}$) that minimizes the sum of the squares of the differences between the observed data points and the values predicted by the model such that $y(t, p)$ satisfies the system of ODEs. Here, $n$ is the number of observation points, $k$ is the number of equations in the ODE system, $\hat{y}_{j}(t_{i})$ denotes the observed data, $y_{j}(t_{i}, p)$ denotes the model-predicted values, and $\sigma^2_{j}$ is the variance of the observed data for component $j$. We then pass this loss function to the Differential Evolution (DE) algorithm \cite{de}, a global optimization method available in the DEoptim library in R \cite{deoptim}. Within the loss function, the system of differential equations is solved using the Livermore Solver for Ordinary Differential Equations (LSODA), which features automatic method switching \cite{liver}. This allows the solver to automatically determine whether the system is stiff and choose the appropriate integration method, removing the need for the user to make that determination.

\begin{figure}[ht]
\centering
\includegraphics[width=18cm]{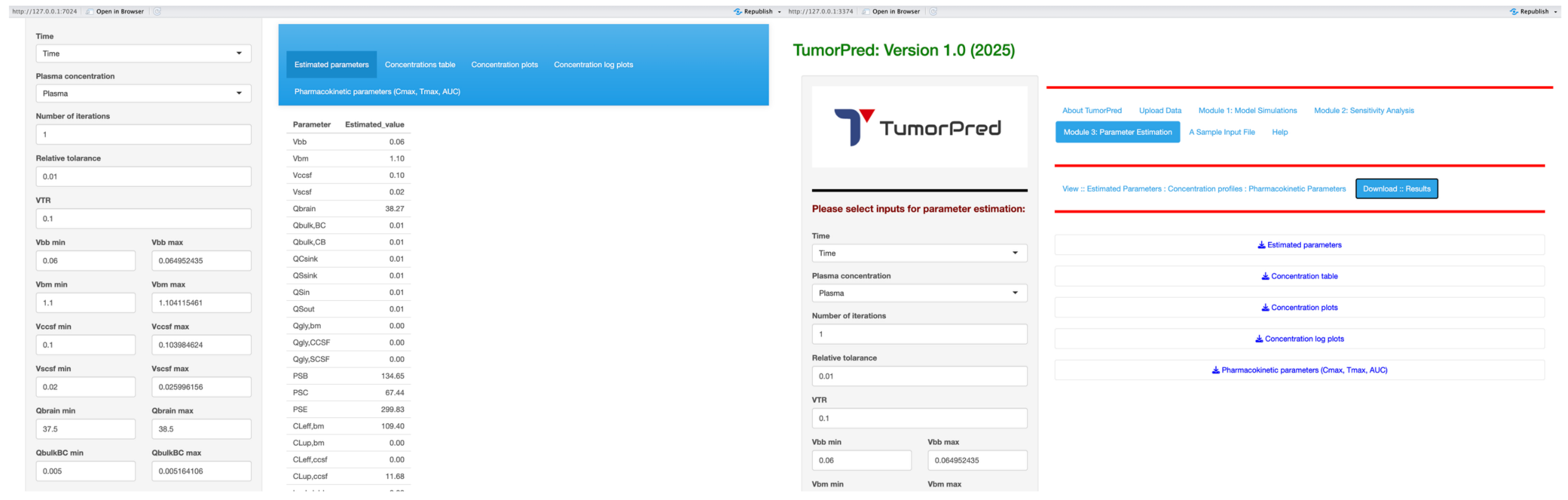}
\begin{tiny}
\caption{A list of TumorPred outputs after a successful parameter estimation.}\label{inputk}
\end{tiny}
\end{figure}

TumorPred is primarily designed to estimate all 27 parameters of the 4-CNS brain model. However, estimating all 27 parameters simultaneously can be computationally expensive. The solution may sometimes fail to converge due to several factors: the high stiffness of the system of ODEs, the high sensitivity of drug-specific and patient-specific parameters, numerical errors from finite computer precision, or a lack of high-quality data. Furthermore, the computational time depends on the number of iterations, the relative tolerance, and the value-to-reach (VTR) termination threshold for the loss function. Values for some drug- and patient-specific parameters can be obtained from experimental data or literature. In such cases, these values can be entered as the minimum and maximum bounds for the parameters via the GUI, which reduces the algorithm's computational burden. Following a successful simulation, the user can view and download a list of the estimated model parameters, the corresponding concentration-time profiles, log concentration-time profiles, and key pharmacokinetic parameters (Cmax, Tmax, and AUC), as shown in Figure~(\ref{inputk}).

\section{Software validation}
For validation purposes, six drug- and patient-specific parameters ($V_{bb}$, $V_{bm}$, $V_{ccsf}$, $V_{scsf}$, $fu_{bb}$, and $\lambda_{ccsf}$) were selected for estimation using the TumorPred application. Table~(\ref{tab:5}) summarizes the parameter estimation results, including the reference values (Ref.), the estimated values (Est.), and the corresponding absolute errors ($ \mathbf{Error} = |\mathbf{Ref} - \mathbf{Est}|$). The results indicate that TumorPred accurately estimated the majority of the selected parameters, demonstrating strong agreement with the reference values.

\begin{figure}[H]
\centering
\includegraphics[width=18cm]{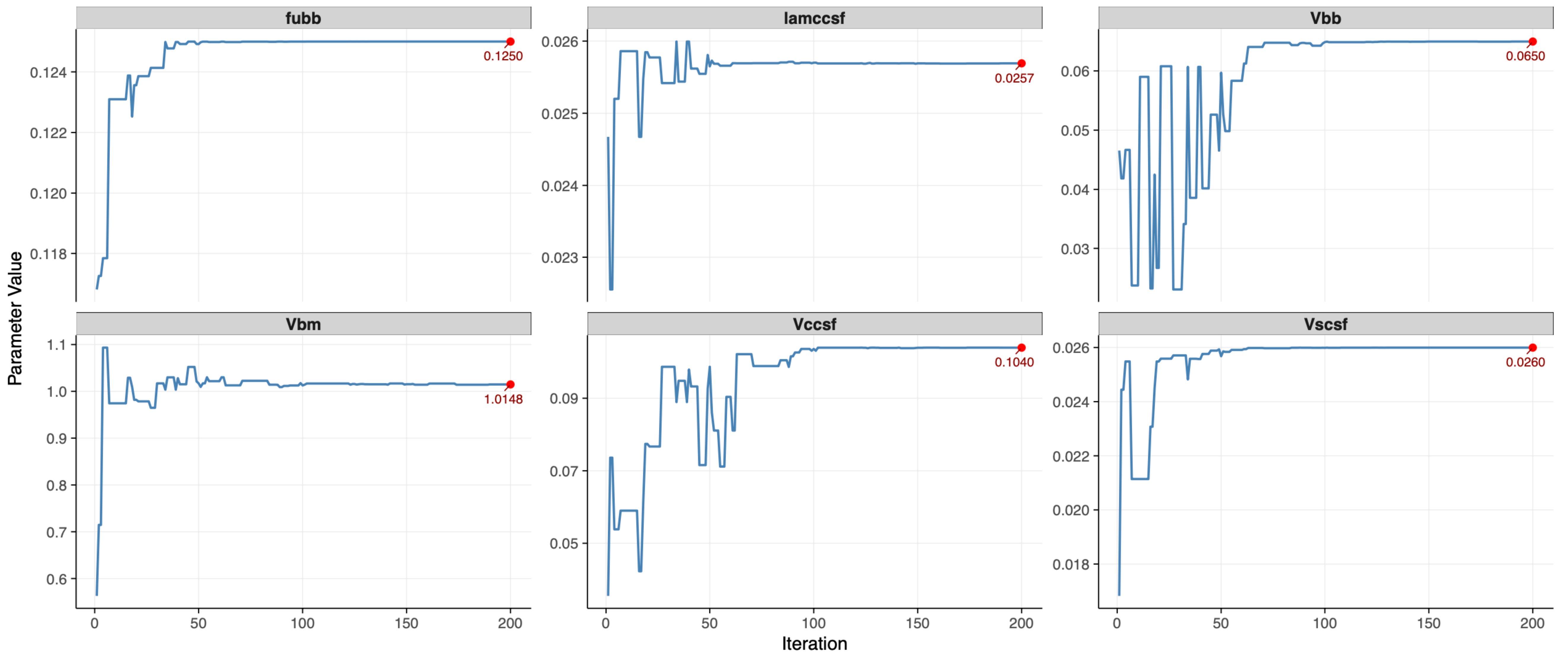}
\begin{tiny}
\caption{Behavior of individual parameter convergence across iterative optimization steps.}\label{pco1}
\end{tiny}
\end{figure}

To further assess the robustness of the estimation procedure, diagnostic plots were generated to illustrate the convergence behavior of the parameters with respect to the number of iterations (Figure~\ref{pco1}). The convergence profiles show progressive stabilization of parameter values, providing evidence that the minimization algorithm implemented in TumorPred successfully identified the global minimum within 200 iterations.

\begin{figure}[ht]
\centering
\includegraphics[width=18cm]{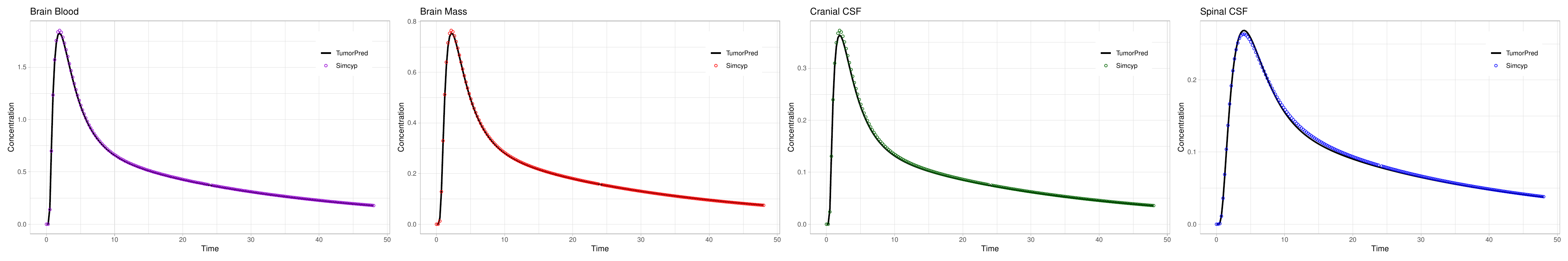}
\begin{tiny}
\caption{Predicted concentration profiles via TumorPred vs reference data obtained from Simcyp Simulator.}\label{pco}
\end{tiny}
\end{figure}

As a secondary level of validation, the estimated parameter values were subsequently used to simulate drug concentration profiles within the TumorPred framework. These simulated profiles were then compared with the reference drug concentration profiles obtained from the Simcyp simulator (Figure~\ref{pco}). The comparison demonstrates close agreement between the estimated and reference concentration profiles.

Collectively, these findings confirm the numerical accuracy and reliability of the TumorPred framework and highlight its potential as a robust computational tool for drug development and optimization of therapeutic strategies in brain cancer treatment.

\begin{table}[htbp]
\centering
\def\arraystretch{1.9}
\caption{Parameter comparison: reference vs. predicted values with absolute errors.} \label{tab:5}
\tabcolsep=12pt
\scriptsize
\begin{tabular}{|c|cccccc|}
\hline
\cline{2-7}
\textbf{Parameter.} & $V_{bb}$ & $V_{bm}$ & $V_{ccsf}$ & $V_{scsf}$ & $fu_{bb}$ & $\lambda_{ccsf}$ \\
\hline
\textbf{Ref.} & \multicolumn{1}{c}{0.064952435} & \multicolumn{1}{c}{1.104115461} & \multicolumn{1}{c}{0.103984624} & \multicolumn{1}{c}{0.025996156} & \multicolumn{1}{c}{0.125000000} & \multicolumn{1}{c|}{0.026000000} \\
\textbf{Est.} & \multicolumn{1}{c}{0.064952243} & \multicolumn{1}{c}{1.014766461} & \multicolumn{1}{c}{0.103953204} & \multicolumn{1}{c}{0.025996112} & \multicolumn{1}{c}{0.124999861} & \multicolumn{1}{c|}{0.025691588} \\
\textbf{Error} & \multicolumn{1}{c}{$1.92\times10^{-7}$} & \multicolumn{1}{c}{$8.93\times10^{-2}$} & \multicolumn{1}{c}{$3.14\times10^{-5}$} & \multicolumn{1}{c}{$4.40\times10^{-8}$} & \multicolumn{1}{c}{$1.39\times10^{-7}$} & \multicolumn{1}{c|}{$3.08\times10^{-4}$} \\
\hline
\end{tabular}
\end{table}

\section{Impact and future work}
The TumorPred application was initially introduced, validated, and partially utilized in two published papers, \cite{Li2} and \cite{Li3}. In \cite{Li2}, the authors developed a novel 9-CNS brain model to measure drug penetration. Subsequently, in \cite{Li3}, they created an R Shiny web application called SpatialCNS-PBPK to make this 9-CNS brain model accessible to a wider audience. To validate the 9-CNS brain model in both studies, the authors compared results from the TumorPred app (which uses a 4-CNS brain model) with those from the Simcyp Simulator v18 \cite{Jamei} (a commercial software which also uses a 4-CNS brain model); this comparison verified the accuracy of TumorPred's calculations. The authors then packaged the same computational mechanism from the TumorPred application into the SpatialCNS-PBPK application to evaluate it through concentration-time profiles, sensitivity analysis, and pharmacokinetic parameter (AUC, Cmax, Tmax) calculation. However, TumorPred has an additional feature for estimating drug- and patient-specific parameters, which was beyond the scope of \cite{Li2} and \cite{Li3}. TumorPred was later fully utilized in the paper \cite{nel} to estimate drug- and patient-specific parameters, drug concentration-time profiles, pharmacokinetic parameters (AUC, Cmax, Tmax), and sensitivity analysis.

As demonstrated in Table 1 and Figure 4, TumorPred is capable of accurately estimating model parameters. For simplicity, this study focused on estimating six key parameters; however, the framework is designed to estimate all drug- and system-specific parameters simultaneously. Although full parameter estimation enhances modeling flexibility, it is computationally intensive. Therefore, future work will focus on optimizing the source code to improve computational efficiency.

In physiologically based pharmacokinetic (PBPK) compartmental brain modeling, system- and drug-specific parameters are often available in the literature or can be obtained through experimental studies. Nevertheless, in many practical situations, parameter values may be unknown or difficult to measure. In such cases, TumorPred provides a valuable solution by enabling data-driven parameter estimation. It is also important to note that although abemaciclib was used as the drug of interest in this study, the framework is generalizable and can be applied to simulations involving other brain cancer therapeutics.

Furthermore, while this work implemented a four-compartment PBPK brain model, the underlying methodology can be extended to other PBPK brain modeling structures reported in the literature \cite{Li4}. Future studies will extend the proposed framework to other PBPK brain models and investigate its application in a wider range of pharmacokinetic and oncology applications.

After accurate parameter estimation, TumorPred can successfully predict drug concentration profiles in the human brain, as illustrated in Figure 5. These predicted concentration profiles allow the calculation of key pharmacokinetic parameters such as AUC, Tmax, and Cmax. The quantitative results obtained from these simulations can be used to guide the design of efficient clinical trials, support the selection of appropriate drug candidates, and optimize dosing regimens. Moreover, the sensitivity analysis feature added to TumorPred greatly helps scientists understand how much the drug concentration changes in response to small variations in model parameters as examine in the Figure 2. Ultimately, this framework supports rational drug development and improves treatment strategies for brain cancer. The TumorPred was made  easily and freely accessible to the user thus accessible to a wider audience. 

\vskip 0.3cm
\textbf{CRediT authorship contribution statement}\\
Charuka Wickramasinghe: Writing – original draft, Software, Conceptualization, Formal analysis, Visualization, Validation, Data curation. Nelum Hapuhinna: Writing – review \& editing, Software, Validation, Visualization, Data curation.
\vskip 0.3cm
\textbf{Declaration of Generative AI and AI-assisted technologies in the writing process}\\
During the preparation of this work, the authors used ChatGPT to enhance their English grammar writing. After using this tool/service, the authors reviewed and edited the content as needed and took full responsibility for the content of the publication.
\vskip 0.3cm
\textbf{Declaration of competing interest}\\
The authors declare that they have no known competing interests that could have appeared to influence the work reported in this paper.

\vskip 0.3cm
\textbf{Acknowledgements}\\
We thank the open-source contributors to R shiny and key R packages whose tools power the underlying simulation. This research received no funding. \\

\end{document}


\section{System of differential equations}

The 4-compartment permeability limited brain (4Brain) model consisting of brain blood, brain mass, cranial and spinal cerebrospinal fluid (CSF). The equations (\ref{eqn1}), (\ref{eqn2}), (\ref{eqn3}), and (\ref{eqn4}) represent the rate of change of drug concentration of brain blood ($C_{bb}$), brain mass ($C_{bm}$), cranial CSF ($C_{ccsf}$), and spinal CSF ($C_{scsf}$) respectively where $C_{art}$ denotes the arterial blood concentration(mg/L) which is an exogenous input to the system.
\vspace{0.5cm}

\noindent \textbf{\textit{Brain blood compartment:}} 
\begin{equation} \label{eqn1}
\begin{split}
V_{bb} \frac{dC_{bb}}{dt} & = Q_{brain}(C_{art}-C_{bb}) + PSB(\lambda_{bm}fu_{bm}C_{bm}-\lambda_{bb}fu_{bb}C_{bb})   \\ & +  CLB_{in}fu_{bb}C_{bb} + CLB_{out}fu_{bm}C_{bm} 
\\ & + PSC(\lambda_{ccsf}fu_{ccsf}C_{ccsf}-\lambda_{bb}fu_{bb}C_{bb}) - CLC_{in}fu_{bb}C_{bb} 
\\ & + CLC_{out}fu_{ccsf}C_{ccsf} + Q_{csink}C_{ccsf} 
+ Q_{ssink}C_{scsf}
\end{split}
\end{equation}
\noindent \textbf {\textit{Brain mass compartment:}}
\begin{equation} \label{eqn2}
\begin{split}
\begin{aligned}
 V_{bm} \frac{dC_{bm}}{dt} &= PSB(\lambda_{bb}fu_{bb}C_{bb}-\lambda_{bm}fu_{bm}C_{bm})  +CLB_{in}fu_{bb}C_{bb} \\& - CLB_{out}fu_{bm}C_{bm} - Q_{bulk}C_{bm}\\
&  + PSE(\lambda_{ccsf}fu_{ccsf}C_{ccsf}-\lambda_{bm}fu_{bm}C_{bm})  - CL_{met}C_{bm} 
\end{aligned}
\end{split}
\end{equation}
\noindent \textbf {\textit{Cranial CSF compartment:}}
\begin{equation} \label{eqn3}
\begin{split}
\begin{aligned}
 \qquad V_{ccsf} \frac{dC_{ccsf}}{dt} 
 &= PSC(\lambda_{bb}fu_{bb}C_{bb}-\lambda_{ccsf}fu_{ccsf}C_{ccsf}) 
+CLC_{in}fu_{bb}C_{bb}  \\&  - CLC_{out}fu_{ccsf}C_{ccsf}  - Q_{sout}C_{scsf} 
 \\&  + PSE(\lambda_{bm}fu_{bm}C_{bm}-\lambda_{ccsf}fu_{ccsf}C_{ccsf})  - Q_{sin}C_{ccsf} 
 \\& - Q_{csink}C_{ccsf}
\end{aligned}
\end{split}
\end{equation}
\noindent \textbf {\textit{Spinal CSF compartment:}}
\begin{equation} \label{eqn4}
\begin{split}
\begin{aligned}
\qquad V_{scsf} \frac{dC_{scsf}}{dt} &= Q_{sin}C_{ccsf} - Q_{sout}C_{scsf} - Q_{ssink}C_{scsf} \hspace{8cm}
\end{aligned}
\end{split}
\end{equation}

\newpage
\section{Parameter descriptions and the reference values}

The mean values of the  system specific (patient specific) parameters  were derived from the Simcyp Simulator for a virtual patient population of 100 patients and listed in the Table \ref{system}. These parameters are related to the physiological and biochemical characteristics of the human body. The drug-specific parameters can be found in the literature and from the Simcyp simulator. The parameter values used in this study are listed in Table \ref{drug}. These parameters are intrinsic to the drug itself and dictate how it behaves in the body. These parameter values are considered as reference values to validate our parameter estimation approach. The drug specific parameters are specific to an oral dose of 10 mg of abemaciclib drug, which is a targeted cancer therapy used to treat certain types of cancer.

\begin{table}[H]
\centering
\caption{System specific parameters for abemaciclib in 4-compartment brain model}
\label{system}
\footnotesize
\setlength{\tabcolsep}{8pt}
\renewcommand{\arraystretch}{1.2}
\begin{tabular}{|>{\RaggedRight}p{0.18\linewidth}|>{\RaggedRight}p{0.5\linewidth}|c|}
\hline
\textbf{Parameter} & \textbf{Description (unit)} & \textbf{Value} \\
\hline
$V_{bb}$ & Brain blood volume (L) & 0.064952435 \\
\hline
$V_{bm}$ & Brain mass volume (L) & 1.104115461 \\
\hline
$V_{ccsf}$ & Cranial CSF volume (L) & 0.103984624 \\
\hline
$V_{scsf}$ & Spinal CSF volume (L) & 0.025996156 \\
\hline
$Q_{brain}$ & Blood/CSF flow (L/h) & 38.0 \\
\hline
$Q_{csink}$ & Cranial CSF absorption rate (L/h) & 0.01277633 \\
\hline
$Q_{ssink}$ & Spinal CSF absorption rate (L/h) & 0.007761342 \\
\hline
$Q_{bulkBC}$ & Bulk flow from brain mass to cranial CSF (L/h) & 0.005164106 \\
\hline
$Q_{bulkCB}$ & Bulk flow from cranial CSF to brain mass (L/h) & 0.005164106 \\
\hline
$Q_{sout}$ & CSF flow: spinal to cranial CSF (L/h) & 0.007489995 \\
\hline
$Q_{sin}$ & CSF flow: cranial to spinal CSF (L/h) & 0.015251337 \\
\hline
$PSB$ & Passive permeability-surface area of BBB (L/h) & 135.0 \\
\hline
$PSC$ & Passive permeability-surface area of BBB (L/h) & 67.5 \\
\hline
$PSE$ & Passive permeability-surface area of brain-CSF barrier (L/h) & 300.0 \\
\hline
\end{tabular}
\end{table}

\begin{table}[H]
\centering
\caption{Drug-specific parameters for abemaciclib in 4-compartment brain model}
\label{drug}
\footnotesize
\setlength{\tabcolsep}{12pt} 
\renewcommand{\arraystretch}{1.2}
\begin{tabular}{|>{\RaggedRight}p{0.18\linewidth}|>{\RaggedRight}p{0.5\linewidth}|c|}
\hline
\textbf{Parameter} & \textbf{Description (unit)} & \textbf{Value} \\
\hline
$CLB_{in}$ & Clearance of active uptake transporter on the BBB (L/h) & 0.0 \\
\hline
$CLB_{out}$ & Clearance of active efflux transporter on the BBB (L/h) & 110.0 \\
\hline
$CLC_{in}$ & BCSFB uptake transporter (L/h) & 11.9 \\
\hline
$CLC_{out}$ & BCSFB efflux transporter (L/h) & 0.0 \\
\hline
$CL_{met}$ & Metabolic clearance due to brain enzymes (L/h) & 0.0 \\
\hline
$fu_{bb}$ & Drug unbound fraction in the brain blood & 0.125 \\
\hline
$fu_{bm}$ & Drug unbound fraction in the brain mass & 0.044 \\
\hline
$fu_{ccsf}$ & Drug unbound fraction in the Cranial CSF & 1.0 \\
\hline
$\lambda_{bb}$ & Unionization fraction in the brain blood & 0.033 \\
\hline
$\lambda_{bm}$ & Unionization fraction in the brain parenchyma & 0.017 \\
\hline
$\lambda_{ccsf}$ & Unionization fraction in the cranial CSF & 0.026 \\
\hline
\end{tabular}
\end{table}
\newpage